\documentclass[12pt, a4paper]{amsart}
\usepackage{floatflt,graphicx,subfigure}
\usepackage{version}
\usepackage{a4wide}
\usepackage{ifthen}

\usepackage{amsmath}
\usepackage{amssymb}

\newcommand\proofbox{\ensuremath{\blacksquare}\relax}
\newcommand{\noproofbox}{\ensuremath{\square}\relax}

\newcommand\singleqedd{\hskip10000pt minus 1fil}
\newcommand\esingleqedd{\hskip10000pt minus 1fil \llap{\noproofbox}}

\def\proof{\par\smallskip
             \noindent {\sc Proof. }}
\def\proofof #1 {\par\medskip\noindent {\sc Proof of #1. }}

\def\sketchof #1 {\par\medskip\noindent {\sc Sketch of proof of #1. }}

\def\qed{\rule{0pt}{0pt}\hfill $\blacksquare$ \par\medskip}
\def\qedd{\rule{0pt}{0pt}\hfill $\square$}

\newenvironment{remark}[1][]{%
  \par\noindent \textsc{\ifthenelse{\equal{#1}{}}{Remark. }{#1. }}}{%
   \par\medskip}

\usepackage{hyperref}

\newcommand{\N}{\mathbb{N}}

\newcommand{\C}{\ensuremath{\mathbb{C}}}

\newcommand{\dt}[1]{\ensuremath{{\mathcal{#1}}}}

\newcommand{\re}{\operatorname{Re}}
\newcommand{\im}{\operatorname{Im}}

\newcommand{\cl}[1]{\overline{#1}}

\newcommand{\Ek}{E_{\kappa}}

\newcommand{\extaddress}[1]{\underline{#1}}

\newcommand{\G}{\mathcal{G}}

\newcommand{\bdyit}[2]
             {{\rule{0pt}{0pt}_{\mbox{$\scriptstyle #2$}}^{\mbox{%
                   $\scriptstyle #1$}} }}

\newcommand{\adds}{\extaddress{s}}
\newcommand{\s}{\adds}

\renewcommand{\r}{\extaddress{r}}

\newcommand{\gs}{g_{\adds}}

\newcommand{\W}{\dt{W}}

\newcommand{\ts}{t_{\s}}

\newcommand{\wt}[1]{\widetilde{#1}}

\newcommand{\PRS}[1]{G_{#1}}
\newcommand{\PR}{\PRS{\s}}

\newtheorem{thm}{Theorem}
\newtheorem{prop}{Proposition}
\newtheorem*{conj}{Conjecture}

\newcommand{\picturedir}{.}
\newcommand{\lowres}{_lr}

\title[A Landing Theorem for Exponential Maps]{%
A Landing Theorem  for Periodic Rays \\ of Exponential Maps}

\author{Lasse Rempe}
\address{Mathematisches Seminar der CAU Kiel, Ludewig-Meyn-Str.~4,
24098 Kiel, Germany}
\curraddr{Mathematics Institute, Warwick University, Coventry CV4 7AL,
United Kingdom}
\email{lasse@maths.warwick.ac.uk}


\date{\today}

\begin{document}

 \begin{abstract}
  For the family of exponential maps $z\mapsto
  \exp(z)+\kappa$, we show the following analog of a theorem of
  Douady and Hubbard concerning polynomials.
  Suppose that $g$ is a periodic dynamic ray of an exponential map
  with nonescaping singular value. Then $g$ lands at a
  repelling or parabolic periodic point. We also show that there are
  periodic dynamic rays landing at all
  periodic points of such an exponential map, with the exception of at
  most one periodic orbit.
 \end{abstract}
 \maketitle
 \section{Introduction}
  In polynomial dynamics, \emph{dynamic rays}, which foliate the
  basin of the superattracting fixed point $\infty$, play an important
  role. By a result of Douady and Hubbard, every periodic dynamic ray of a
  polynomial with connected Julia set lands at a repelling or
  parabolic periodic point \cite[Theorem 18.10]{jackdynamics}.
  The combinatorial information given by this fact was a
  key element of their initial study of the Mandelbrot set
  \cite{orsay}, and has since proved to be one of the most fundamental 
  tools in the field.

  In the case of exponential maps 
  $E_{\kappa}(z) = \exp(z) + \kappa$, the point $\infty$ is no
  longer an attracting fixed point but rather an essential
  singularity, and the set of escaping points is no longer open.
  However, it has long been known \cite{dgh} that there exist curves
  of escaping points in the dynamical plane which can be seen as an
  analog of dynamic rays (compare Section \ref{sec:dynamicrays}). 
  The combinatorics of such rays landing
  together has also been used to great advantage in this family (see,
  e.g., \cite{expattracting, boundary,expcombinatorics}).

 Therefore it is a natural and important question 
  whether periodic dynamic rays of exponential maps
  land in general; compare \cite[Section VI.6]{habil}. In this article, we
  give a positive answer to this question.

\begin{thm} \label{thm:landing2}
 Let $\kappa\in\C$ such that the singular value $\kappa$ of $\Ek$ does not
 escape to infinity. Then every periodic dynamic ray of $\Ek$ lands at a
 repelling or parabolic periodic point.
\end{thm}
\begin{remark}
 For maps with escaping singular value, there is an obvious
 exception if the singular value lies on a periodic dynamic ray (and this is
 the only exception). Since these parameters were already treated in
 \cite{expper}, we
 have formulated the theorem for nonescaping parameters to simplify
 its statement.
\end{remark}

The analog of Theorem \ref{thm:landing2} in the polynomial case is
 proved by a hyperbolic contraction argument. In the exponential
 family, this proof can be modified to work in the special case where
 the ray considered
 does not intersect the postsingular set
  \cite{expper,thesis}. 
  This condition is satisfied for all periodic dynamic
  rays in the important cases of
  attracting, parabolic
  and postsingularly finite parameters; see
  \cite{expper}. For other
  conditions under which the landing of periodic rays
  can be proved by an argument in the dynamical plane, compare
  \cite{landingsiegel,accessiblesingularvalue}.

On the other hand, there is a generic subset of the
  bifurcation locus consisting of maps whose
  postsingular set
  is the entire plane 
  \cite[Theorem 5.1.6]{thesis}. In this case, it is
  impossible to prove by
  hyperbolic contraction that a given periodic ray lands,
  unless one can \emph{a priori} ensure that the ray does not
  accumulate on the singular value $\kappa$. However, the accumulation
  behavior of dynamic rays can be very complicated in general, e.g.,
  for a large set of parameters there exist rays whose accumulation sets
  are indecomposable continua
  \cite{nonlanding}. It
  even seems difficult to show dynamically, without any structural
  assumptions, that the accumulation set of a periodic ray cannot be
  the entire plane.

Our proof of Theorem \ref{thm:landing2}
 manages to circumvent these problems by working in the parameter
 plane, using
 the ``$\lambda$-lemma'' from the
 theory of holomorphic motions \cite{mss}
 and a theorem of Schleicher concerning the landing of periodic
 \emph{parameter} rays \cite{habil}. 
 It is somewhat surprising that the landing of
 dynamic rays is proved using landing properties of parameter rays.
 Unfortunately, this method will not carry
 over to higher-dimensional parameter
  spaces of entire transcendental functions.

 We also prove a partial converse to Theorem \ref{thm:landing2}. In
  the case of a polynomial with connected Julia set, it is known that
  every repelling or parabolic periodic point is the landing point of
  a dynamic ray \cite[Theorem 18.11]{jackdynamics}. The analog of this
  theorem has previously been shown
  for attracting, parabolic, escaping and postsingularly finite
  exponential maps \cite{expper}. (In the escaping case,
  there will again be exceptions when the singular value lies on a periodic
  ray.)

\begin{thm} \label{thm:perpoints}
 Let $\kappa\in\C$
  such that the singular value $\kappa$ of $\Ek$ does not escape. Then, with
  the exception of at most one periodic orbit, every
  periodic point of $\Ek$ is the landing point of a periodic dynamic ray.
\end{thm}

 An irrationally indifferent periodic point can never be the landing
 point of a periodic dynamic ray by the Snail Lemma \cite[Lemma
 16.2]{jackdynamics}. Thus our theorem implies that, in the case
 where $\kappa$ is a Siegel or Cremer parameter,
 all repelling periodic points are landing points.
 We believe that, as for polynomials,
 nonrepelling orbits are the only
 exceptions which can actually occur in Theorem \ref{thm:perpoints}. 
 This conjecture can be phrased as a question about
 parameter space (see Section \ref{sec:perpoints}) and would imply
 that every indifferent parameter $\kappa$ is the landing point of a
 parameter ray.

\subsection*{Organization of the article}
 Naturally, we will rely on a fair number of results now established in
  the theory of exponential maps. We have attempted to collect these
  in Sections \ref{sec:dynamicrays} and \ref{sec:wakes}, so that the
  proofs of Theorems \ref{thm:landing2} and \ref{thm:perpoints} will only
  rely on the propositions stated there.

\subsection*{Acknowledgments} I would like to thank Walter Bergweiler and
  Dierk Schleicher, as well as the audiences of the seminars at
  Nottingham, Open University and Warwick for interesting discussions on
  this work. Much appreciation is owed to Misha Lyubich, John Milnor and the
  Institute of Mathematical Sciences at Stony Brook for continued
  support and encouragement. I also wish to thank
  the University of Warwick, where part of this work was
  conducted while holding a Marie Curie Studentship, for its hospitality.

\subsection*{Notation and preliminaries.} 
  Throughout this article, $\C$ will denote the complex plane,
 and all closures will be taken in $\C$ (rather than on the Riemann
 sphere). The set of \emph{escaping points} of an exponential map
 $\Ek$ is denoted by
  \[ I(\Ek) := \{z\in\C: |\Ek^n(z)|\to\infty\}. \]

 By \cite{alexmisha}, an exponential map has at most one periodic orbit
  which is not repelling. Therefore, we call a
  parameter $\kappa$ 
  \emph{attracting}, resp.\ \emph{parabolic}, if $\Ek$ has an
  attracting, resp.\ parabolic (i.e., rationally indifferent), periodic orbit.
  As a model of exponential growth, we will fix the function 
  $F:[0,\infty)\to [0,\infty), F(t) = \exp(t) - 1$.

 The completion 
  of a proof is marked by the symbol $\proofbox$. Results which are
  cited without proof are concluded by $\noproofbox$.

\section{Dynamic and Parameter Rays of Exponential Maps} 
          \label{sec:dynamicrays}

\subsection{Dynamic rays}
We call a
 sequence $\s=s_1 s_2 s_3 \dots$ of integers an \emph{external
 address}.
 We say that a point $z\in\C$ \emph{has external address $\s$} if
 \[ \im \bigl( \Ek^{n-1}(z)\bigr)\in \bigl((2s_n -1)\pi,(2s_n+1)\pi\bigr) \] 
 for every
 $n\geq 1$. 
The shift map on external addresses is denoted by
 $\sigma$; i.e. $\sigma(s_1 s_2 s_3 \dots )=s_2 s_3 s_4 \dots$\ .
 For simplicity, we will restrict to \emph{bounded}
 external addresses; i.e.\ those for which the sequence $(s_k)$ 
 of entries is bounded. For the results in the general case, see
 \cite{expescaping} or \cite{topescapingnew}.

 \begin{prop}[\cite{expescaping}] \label{prop:dynamicrays}
  Let $\kappa\in\C$, and let $\s$ be a bounded external address. 
   If $T>0$ is large enough, then there exists a unique
   curve
   $g:[T,\infty)\to I(\Ek)$ such that
   \begin{enumerate}
    \item $g(t)$ has external address
     $\s$ for all $t\geq T$, and \label{item:extaddress}
    \item $\re\bigl( \Ek^n(g(t))\bigr)=
          F^n(t)+o(1)$ as $n$ or $t$ tend to $\infty$.
          \label{item:asymptotics}
  \end{enumerate}
  This curve has a unique maximal extension
   $g^{\kappa}_{\s}:(t^{\kappa}_{\s},\infty)\to I(\Ek)$ 
   satisfying (\ref{item:asymptotics}). This extension 
   is injective and will be called the
   \emph{dynamic ray of $\Ek$ at address $\s$}.

  The family of dynamic rays has the following additional properties:
   \begin{itemize}
    \item $\Ek(g^{\kappa}_{\s})\subset
               g^{\kappa}_{\sigma(\s)}$ for all $\s$;
    \item for every $\s$, $t^{\kappa}_{\s}$
           depends upper semicontinuously on $\kappa$; 
    \item for every $t>0$, 
      $g^{\kappa}_{\s}(t)$ depends
      analytically on $\kappa$ (where defined);
    \item
      $t_{\s}^{\kappa}\neq 0$ if and only if there exists
     $n\geq 1$ with 
     $\kappa\in g^{\kappa}_{\sigma^n(\s)}$ (in which case we call
     $g^{\kappa}_{\s}$ \emph{broken}).\qedd
   \end{itemize}
 \end{prop} 
\begin{remark}
  In the case of unbounded addresses, the number $0$ in the last
  item is replaced by a number $t_{\s}\geq 0$ which depends only
  on $\s$. 
\end{remark}

 As usual, we say that a dynamic ray $\gs^{\kappa}$ \emph{lands} if
  $t^{\kappa}_{\s}=0$ and $\lim_{t\to 0} g_{\s}^{\kappa}(t)$ exists. 
  In our proofs of 
  Theorems \ref{thm:landing2} and \ref{thm:perpoints}, we will use the
  fact, mentioned in the introduction, that these theorems
  are known to hold for attracting and parabolic parameters.

 \begin{prop}[{\cite[Theorem 3.2 and Theorem 5.4]{expper}}]
  \label{prop:landing1}
  Let $\kappa$ be an attracting or parabolic parameter. Then every
   periodic dynamic ray of $\Ek$ lands at a repelling or parabolic
   periodic point. Conversely, every repelling or parabolic periodic point
   of $\Ek$ is the landing point of a periodic dynamic ray. \qedd
 \end{prop}
 \begin{remark}
   In fact, for these parameters even more is true:
  every dynamic ray lands, and every point in the Julia set is either on
  a ray or the landing point of a ray \cite{accessible, topescapingnew}. 
 \end{remark}

 \subsection*{Parameter rays}
 By \cite[Theorem~II.8.1]{habil} (which was recently generalised
   to the case of unbounded addresses; see
   \cite{markus,markusdierk,markuslassedierk}),
   the set
   \[ \PR := \bigl\{\kappa\in\C: 
               \kappa\in g^{\kappa}_{\s}\bigl((t_{\s}^{\kappa},\infty)\bigr)
             \bigr\} \]
   is itself a continuous (and differentiable) curve to $\infty$
   for every bounded external address $\s$. This curve is
   called the \emph{parameter ray} at address $\s$. We shall require
  the following result of 
  Schleicher, stating 
   that every periodic parameter ray lands. 
   (The \emph{ray period} of a
   parabolic parameter $\kappa$ is the period of the repelling
   petals of
   $E_{\kappa}$ at the parabolic orbit; in particular, it is a multiple
   of the period of this orbit.)

 \begin{prop}[\protect{\cite[Theorem V.7.2]{habil}}] \label{prop:parraysland}
  Suppose that $\s$ is a periodic address of period $n$.
  Then $\overline{\PR}\setminus \PR$ consists of a single
  parabolic parameter of ray period $n$. 

 Conversely,
  for every parabolic parameter $\kappa$ of ray period $n$, there exists
  a periodic address of period $n$ such that
  $\kappa\in\overline{\PR}$. 
 \end{prop}
 Schleicher's proof of this proposition consists of two
  steps. First it is shown that every parabolic parameter is the
  landing point of one or two parameter rays (whose addresses depend on the
  combinatorics of the parameter $\kappa$). This result is
  so far unpublished outside of \cite{habil}. It is then shown 
  (independently from the previous step) that,
  for every periodic address $\s$, there exists a
  parabolic parameter with the required combinatorics. A proof
  of this can be found in 
  \cite[Theorem 8.5]{boundary}. 
  \qedd

\subsection*{Vertical order}
 Finally, let us note that the vertical order of dynamic rays
  coincides with the lexicographic order $<$ on external addresses. More
  precisely, let $\s\neq\r$ be bounded external addresses, and let
  $t_0>\ts^{\kappa}$. If $R>0$ is large enough, then the set
  \[ \{z:\re z > R\}\setminus
         \gs^{\kappa}\bigl([t_0,\infty)\bigr) \]
  has two unbounded components $U^+$ and $U^-$, with unbounded
  positive, respectively negative, imaginary parts.
  The ray $g_{\r}$ tends to
  $\infty$ in $U^-$ if and only if $\r<\s$. (This fact follows easily
  from the definition of dynamic rays;
  see \cite[Lemma 3.9]{markusdierk} or
  \cite[Lemma 3.7.1]{thesis} for details.)

 This statement can be transferred to the parameter plane
  \cite[Proposition 3.10]{markusdierk}. That is,
  if $\s^1$ and $\s^2$ are bounded external addresses, then
  $\PRS{\s^1}$ is below $\PRS{\s^2}$ in the above sense if and only if
  $\s^1<\s^2$.

\section{Landing of Periodic Dynamic Rays}  \label{sec:landing2}
 The key idea in our proof of Theorem \ref{thm:landing2}
 is that dynamic rays move 
 holomorphically with respect to the parameter (where defined). By the
 ``$\lambda$-lemma'' from \cite{mss}, the same is true for their closure. 
 For completeness, we include the short proof of this fact
 (which is analogous to the argument for the general case in \cite{mss}).

\begin{thm} \label{thm:lambdalemma}
 Let $\s$ be a bounded external address.
 Let $W$ be a connected component of the set
  \[ \C\setminus \cl{\bigcup_{n\geq 1} \PRS{\sigma^n(\s)}}. \]
 Suppose furthermore that for some
 $\kappa_0\in W$, $g_{\s}^{\kappa_0}$ lands. 
 Then $g_{\s}^{\kappa}$ lands 
 for every $\kappa\in W$, and the landing
 point depends holomorphically on $\kappa$.
\end{thm}
\begin{remark}
  The theorem (and its proof) remains true for arbitrary 
 external addresses.
\end{remark}
\proof Let $z_0$ be the landing point of $\gs^{\kappa_0}$.
 By definition of $W$, the ray $\gs^{\kappa}$ is unbroken 
 (i.e., $\ts^{\kappa}=0$)
 for every $\kappa\in W$.
 We shall use the notation $h_t(\kappa) := g_{\s}^{\kappa}(t)$.
 The family $(h_{t})_{t\in(0,1)}$ of
 holomorphic functions is normal by Montel's theorem, as it omits the
 holomorphic functions $\kappa\mapsto \kappa$ and
 $\kappa\mapsto h_1(\kappa)$. Let $h_0$ be a
 limit function of this family as $t\to 0$.
 Note that $h_0(\kappa_0)=z_0$ and thus, by Hurwitz's theorem,
  $h_0(\kappa)\neq h_t(\kappa)$ for all $\kappa\in W$ and
   $t>0$. 

  Now if $\wt{h_0}$ is some other limit function of $h_t$ as $t\to 0$,
  then $\wt{h_0}(\kappa_0)=z_0=h_0(\kappa_0)$. Since the $h_t$ omit the
  holomorphic function $h_0$, we can apply Hurwitz's theorem
  once more to see that $\wt{h_0}=h_0$. 
  In other words, $h_t\to h_0$ as $t\to 0$. In particular,
  $g_{\s}^{\kappa}(t)\to h_0(\kappa)$ for every $\kappa\in W$. \qed

\proofof{Theorem \ref{thm:landing2}} 
 Let $\s$ be any periodic address of period $n$, and let
  $\kappa\in\C$ be a parameter with $\kappa\notin I(E_{\kappa})$. 
  Set
  \[ K := \bigcup_{1\leq j\leq n} \cl{\PRS{\sigma^j(\s)}}. \]
  By Propositions \ref{prop:parraysland} and
    \ref{prop:landing1}, we may suppose that $\kappa\not\in K$.

 Let $W$ be the
  connected component of $\C\setminus K$ containing $\kappa_0$.
  It is a standard fact
  that $W$ contains some attracting parameter $\kappa_0$.
  (This follows from Proposition~\ref{prop:wake} below, but can be proved in
  a much more elementary way: by \cite[Theorem 3.5]{expattracting},
  there is a hyperbolic component of parameter space
  (compare Section \ref{sec:wakes}) tending to $\infty$ between any two
  parameter rays.) The ray $g_{\s}^{\kappa_0}$ lands by
  Proposition \ref{prop:landing1}, and by
  Theorem~\ref{thm:lambdalemma}, the
  ray $g_{\s}^{\kappa}$ also lands.
  By the Snail Lemma \cite[Lemma 16.2]{jackdynamics},
  its landing point is repelling or parabolic. \qed

\section{Wakes of Hyperbolic Components} \label{sec:wakes}
 The article \cite{expcombinatorics} gives a complete description of
  the combinatorial structure of \emph{wakes} in exponential parameter space.
  Proposition \ref{prop:parraysland} 
  and Theorem \ref{thm:lambdalemma} can be used
  to translate these results into statements about the regions in which
  certain periodic rays land together, which will
  be needed in the proof of Theorem \ref{thm:perpoints}. 
  These results, which are all analogs of
  well-established results for the Mandelbrot set
  (see, e.g., \cite{jackrays}), will be given in this section. 
  Their proofs necessarily involve citing several combinatorial results
  from \cite{expcombinatorics}. 

\begin{figure}%
 \includegraphics[width=\textwidth]{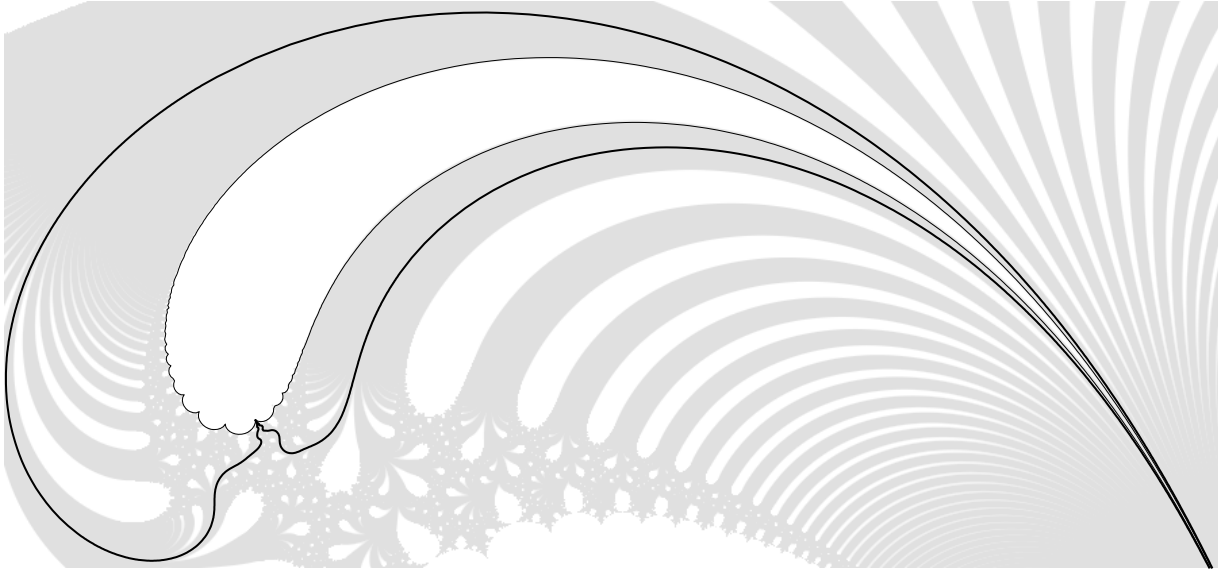}%
 \caption{\label{fig:wake}The wake of a period $5$ hyperbolic component.}
\end{figure}

A \emph{hyperbolic component} $W$ is a maximal connected subset of parameter
space in which all maps have an attracting periodic orbit, the
period of which is necessarily constant on $W$. It is well-known
\cite{bakerexp} that
there is a unique hyperbolic component of period one, which contains the
left half-plane $\{\re \kappa < -1\}$.

\begin{prop} \label{prop:wake}
 Let $W$ be a hyperbolic component of period $n>1$. Then there exists a
 unique parabolic parameter $\kappa_0\in\partial W$ of ray period $n$
 (called the \emph{root} of $W$) 
 which is the common landing point of
 two periodic parameter rays $\PRS{\s^1}$ and $\PRS{\s^2}$
 of period $n$. These rays are called the \emph{characteristic rays} of
 $W$, and they separate $W$ from the unique hyperbolic component of
 period one.
 The connected component of \[
   \C\setminus \bigl(\PRS{\s^1}\cup\{\kappa_0\}\cup \PRS{\s^2}\bigr) \]
   which contains $W$ is called the \emph{wake} of $W$ and denoted by
   $\W(W)$. Any two periodic parameter rays which land at a common
   point are the characteristic rays of some hyperbolic component. \qedd
\end{prop}
\begin{remark}
  The unique period one component does not have
   characteristic rays; by definition, its wake is all of $\C$.
\end{remark}
For a proof of Proposition
 \ref{prop:wake}, compare
 Definition 8.6 and the preceding discussion in
 \cite{boundary}.
 The wake of a hyperbolic component is important to us because
 configurations of periodic dynamic rays landing together are
 preserved throughout it. 
 In order to make this precise, suppose that
 $z_0$ is a periodic point of period $n$ for an exponential map
 $E_{\kappa_0}$ and $U$ is an open neighborhood of $\kappa_0$ in parameter 
 space. By an
 \emph{analytic extension of $z_0$ on $U$} we mean an
 analytic function $z:U\to\C$ with 
 $z(\kappa_0)=z_0$ and $\Ek^n(z(\kappa))=z(\kappa)$ for all $\kappa\in U$.
 By the implicit function theorem and the identity principle,
 $z_0$ has at most one analytic continuation on $U$ provided that
 $(\Ek^n)'(z_0)\neq 1$. Furthermore, if $U$ is simply connected and
 contains no parabolics whose ray period divides $n$, then $z_0$ does have
 an analytic extension to $U$ by the monodromy theorem. 

\begin{prop} \label{prop:landinginwake}
 Let $W$ be a hyperbolic component and $\kappa_0\in W$. Suppose that
 $E_{\kappa_0}$ has a periodic point $z_0$ which is the landing point of
 two different periodic
 dynamic rays $g_{\s^1}^{\kappa_0}$ and $g_{\s^2}^{\kappa_0}$.
 Then $z_0$ has an analytic extension 
 $z$ to $\W(W)$, and for all $\kappa\in\W(W)$, 
  $g_{\s^1}^{\kappa}$ and $g_{\s^2}^{\kappa}$ land at $z(\kappa)$.
\end{prop}
\proof  By 
 \cite[Lemma 5.2]{expattracting}, there exist 
 two periodic addresses $\r^-<\r^+$ such that 
 \begin{enumerate}
  \item
    $g_{\r^-}^{\kappa}$
    and $g^{\kappa}_{\r^+}$ land
    at a common point belonging to the orbit of $z_0$, and 
  \item
    for all $k\geq 0$ and each
     $j\in\{1,2\}$, $\sigma^k(\s^j)\notin (\r^-,\r^+)$.
    \label{item:characteristicproperty}
 \end{enumerate}

 By \cite[Proposition 7.4 and Lemma 3.9]{expcombinatorics}, 
  $\r^-$ and $\r^+$ are the
  characteristic addresses of some hyperbolic component 
  $W'$ with $W\subset \W(W')$. By
  (\ref{item:characteristicproperty}) and the fact that the vertical order of
  parameter rays coincides with the lexicographic order on external
  addresses,
  \[ \W(W')\cap \Biggl(\bigcup_{k\geq 0,j\in\{1,2\}} \G_{\sigma^k(\s^j)}\Biggr)
       = \emptyset. \]
 By Theorem \ref{thm:lambdalemma}, for each $j\in\{1,2\}$ and
  every $\kappa\in \W(W)\subset \W(W')$, the 
  ray $g_{\s^j}^{\kappa}$ lands at a periodic point $z^j(\kappa)$,
  which depends holomorphically on
  $\kappa$. The maps $z^1$ and $z^2$ are both analytic extensions of the
  periodic point $z_0$, and must therefore be equal, as required. \qed

\begin{figure}%
 \includegraphics[width=\textwidth]{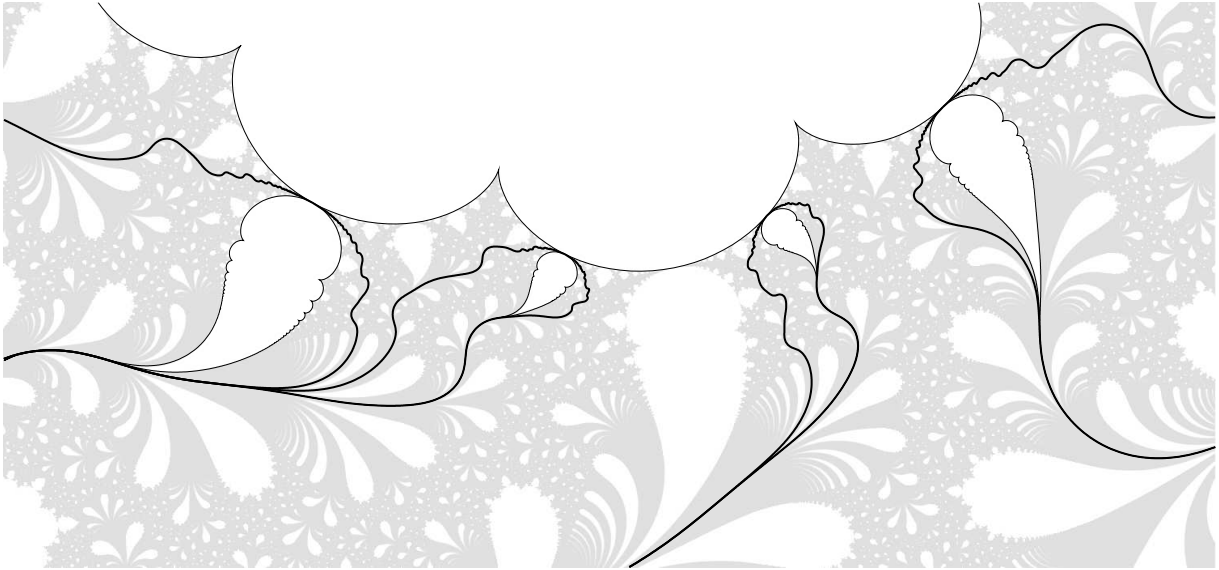}%
 \caption{\label{fig:subwakes} Some subwakes of the period $5$
 component from Figure \ref{fig:wake}.}
\end{figure}

Suppose that $W$ is a hyperbolic component
   and $V$ is a hyperbolic component of
larger period than $W$ such that $\partial V\cap \partial W\neq
\emptyset$.
Then $V$ is called a \emph{child component} of $W$, and
$\W(V)$ is called a \emph{subwake} of $W$ (compare Figure \ref{fig:subwakes}).

\begin{prop} \label{prop:subwake}
 Let $W$ be a hyperbolic component and let $V$ be a child component of $W$.
  Then $\partial W\cap \partial V$ consists only of the root $\kappa_1$ of $V$.

 Furthermore,
  let $\kappa_0\in W$, let
  $z_0$ be a point on the attracting orbit of $E_{\kappa_0}$ and let
  $U$ be an open neighborhood of $W\cup \{\kappa_1\}$ 
  such that $U\cap \W(V)$ is connected and $z_0$ has
  an analytic extension $z$ to $U$. Then $z$ extends analytically to
  $U\cup \W(V)$ and
    $z(\kappa)$ is the landing point of at least two periodic dynamic rays
  for every $\kappa\in \W(V)$. 

 If $W'\neq W$ is a hyperbolic component with $W'\subset\W(W)$, then 
 $W'$ is contained in some subwake of $W$. 
\end{prop} 
\proof According to
 \cite[Proposition 8.1]{boundary}, $\partial W\cap \partial V$ contains
  only parabolic parameters. 
  By \cite[Proposition 5.1 and Theorem 5.3]{expcombinatorics},
  no parabolic parameter except
  $\kappa_1$ can belong
  to this intersection, and
  there are at
  least two periodic rays
  landing at each point of the parabolic orbit of
  $E_{\kappa_1}$. The point
  $z(\kappa)$ belongs to this parabolic orbit when 
  $\kappa=\kappa_1$, and
  becomes repelling when $\kappa$ is perturbed into $V$.

  By \cite[Proposition 5.2]{expcombinatorics}, the
  point $z(\kappa)$ is still
  the landing point of at least two periodic rays after such a perturbation,
  and the second statement now follows from
   Proposition \ref{prop:landinginwake}.
   The final claim is
   \cite[Corollary 6.9]{expcombinatorics}. \qed

Finally, we shall require the following bound for
 parameters on parameter rays with large combinatorics 
 \cite[Corollary 4.8]{topescapingnew}.
 (Also compare
 \cite[Lemma 5.10]{boundary}.)

\begin{prop} \label{prop:parameterraybound}
 Let $\s$ be an external address of period $n$ such
 that
  $M:=\max |s_k| > F^n(6)$. Then
  \[\singleqedd \PR \subset \left\{\kappa:|\kappa| >
                                     \frac{1}{5} F^{-(n-1)}(2\pi
    M) \right\}. \esingleqedd\]
\end{prop}

\section{Proof of Theorem \ref{thm:perpoints}} \label{sec:perpoints}

 \begin{thm} \label{thm:characterization}
  Let $E_{\kappa_0}$ be an exponential map for which the singular value does
  not escape. Suppose that $E_{\kappa_0}$ has a periodic point $z_0$ of period $n$
  which is not the landing point of a periodic dynamic ray. Then
  there exists a hyperbolic component $W$ of period $n$ such that
  $\kappa$ lies in
  $\W(W)$, but not in any subwake of
  $W$.

  Furthermore, all period $n$ points of $E_{\kappa_0}$ which do not lie on
  the orbit of $z_0$ are landing points of periodic dynamic rays.
 \end{thm}
 \proof By Proposition
  \ref{prop:landing1}, 
  we can suppose that $\kappa_0$ is not a parabolic parameter.
  Consider the set $A$ consisting of all parameter rays of
  period dividing $n$ together with the
  landing points of these rays; by Proposition
  \ref{prop:parraysland}, the latter points are exactly the 
  parabolic parameters
  whose ray period divides $n$. We claim that $A$ is closed. 
  Indeed, for every $M\in\N$ the finite union
   \[ A_M := \bigcup\ \bigl\{\,\cl{\PR}\,:\text{$\s$ is periodic of period 
                           dividing $n$ and
                            $|s_j|\leq M$ for all $j$}\bigr\} \]
   is closed. By Proposition \ref{prop:parameterraybound}, the sequence
    of sets
    $(A_M\setminus A_{M-1})$ has no finite accumulation point as
   $M\to\infty$. Thus $A = \bigcup A_M$ is closed as well.

   Now let $U$ be the component of $\C\setminus A$ which contains
   $\kappa_0$. By Proposition \ref{prop:wake}, $U$
   contains a unique hyperbolic component $W$ of period
   $k\vert n$, and this component satisfies
   $U\subset \W(W)$. Since $A\cup\{\infty\}$ is a connected subset of
   the Riemann sphere, $U$ is simply connected.
   Thus any period $n$ point $z$ of $E_{\kappa_0}$
   has an analytic extension
   $z:U\to\C$.
   Let $\kappa\in W$. We claim that if $z$ is not the landing
   point of a periodic dynamic 
   ray of $E_{\kappa_0}$, then $z(\kappa)$ is on the
   attracting orbit of $\Ek$ (and in particular, $k=n$). 
   Indeed, suppose that $z(\kappa)$ is not on the attracting orbit of
   $\Ek$;
   then it is a repelling
   periodic point of $\Ek$. By Proposition
   \ref{prop:landing1}, there is
   a periodic dynamic
   ray $\gs^{\kappa}$ landing at $z(\kappa)$. If more than one
   periodic ray lands at $z(\kappa)$,
   then by
   Proposition \ref{prop:landinginwake}, $\gs^{\kappa_0}$ lands
   at
   $z$. Otherwise, 
   $\s$ has period
   $n$. Because $U$ does not intersect any parameter rays of period
   $n$,
   Theorem \ref{thm:lambdalemma} implies that, again, 
   $\gs^{\kappa_0}$ lands
   at $z$.

 So the analytic continuation of
  $z_0$ throughout $U$ becomes attracting in $W$, 
  and every period $n$ point which is not on the orbit of $z_0$
  is the landing point of a periodic dynamic ray of $E_{\kappa_0}$. 
  If $V$ is a child component of $W$, then clearly $\W(V)\cap U$
  is connected. Thus $\kappa\notin \W(V)$ by Proposition
  \ref{prop:subwake}. \qed

 \proofof{Theorem \ref{thm:perpoints}} 
  Let $z$ and $w$ be periodic points of $\Ek$ which lie on
  two different periodic orbits and let
  $n$ and $m$ be their respective periods.
  Suppose by contradiction that
  neither point is the landing point of a periodic dynamic
  ray of $\Ek$. By the previous theorem, $n\neq m$ and
  there exist period $n$
  and $m$ hyperbolic components $W$ and $V$ with $\kappa\in \W(W) \cap
  \W(V)$ so that $\kappa$ is not contained in any subwake of $W$ or
  $V$. 
  Because the two wakes are not disjoint, one is contained in the
  other; suppose without loss of generality that $\W(W)\subset
  \W(V)$. By Proposition \ref{prop:subwake}, $\W(W)$ lies in a
  subwake of $V$, which is a contradiction. \qed

 By Theorem \ref{thm:characterization}, 
 the conjecture that every repelling periodic point of an exponential
 map with nonescaping singular value is the landing point of a
 periodic dynamic ray is equivalent to the following
 conjecture about
 parameter space. (For a discussion,
 see \cite[Section 7.2]{thesis}.)

 \begin{conj}
  Let $W$ be a hyperbolic component. Then $\W(W)$ is the union
  of $\cl{W}$, the subwakes of $W$, and certain parameter rays.
 \end{conj}

\bibliographystyle{hamsplain}
\bibliography{C:/Latex/Biblio/biblio}

\providecommand{\href}[2]{#2}\input{cyracc.def} \def\j{{\u i}} \def\J{{\u I}}
  \newfont{\cyrit}{wncyi10 at 12pt}\def\cprime{$'$}
\providecommand{\bysame}{\leavevmode\hbox to3em{\hrulefill}\thinspace}
\begin{thebibliography}{BDG}

\bibitem[BR]{bakerexp}
I.~No\"el Baker and Philip~J. Rippon, \emph{Iteration of exponential
  functions Ann},. Acad. Sci. Fenn. Ser. A I Math. \textbf{9} (1984), 49--77.

\bibitem[BDD]{accessible}
Ranjit Bhattacharjee, Robert~L. Devaney, R.~E.~Lee Deville, Kre{\v{s}}imir
  Josi{\'c}, and Monica Moreno-Rocha,
  \emph{\href{http://math.bu.edu/people/bob/papers/accessible.ps}{Accessible
  points in the {J}ulia sets of stable exponentials}}, Discrete Contin. Dyn.
  Syst. Ser. B \textbf{1} (2001), no.~3, 299--318.

\bibitem[DG]{devgoldberg}
Robert~L. Devaney and Lisa~R. Goldberg, \emph{Uniformization of attracting
  basins for exponential maps}, Duke Math. J. \textbf{55} (1987), no.~2,
  253--266.

\bibitem[DGH]{dgh}
Robert~L. Devaney, Lisa~R. Goldberg, and John~H. Hubbard, \emph{A dynamical
  approximation to the exponential map by polynomials}, Preprint, 1986.

\bibitem[DH]{orsay}
Adrien Douady and John Hubbard, \emph{Etude dynamique des polyn{\^o}mes
  complexes}, Pr{\'e}publications math{\'e}mathiques d'Orsay (1984 / 1985),
  no.~2/4.

\bibitem[EL]{alexmisha}
Alexandre~{\`E}. Eremenko and Mikhail~Yu. Lyubich, \emph{Dynamical properties
  of some classes of entire functions}, Ann. Inst. Fourier (Grenoble)
  \textbf{42} (1992), no.~4, 989--1020.

\bibitem[F]{markus}
Markus F{\"o}rster, \emph{Parameter rays for the exponential family},
  Diplomarbeit, Techn. Univ. M{\"u}nchen, 2003, Available as
  \href{http://www.math.sunysb.edu/cgi-bin/thesis.pl?thesis03-3}{Thesis
  2003-03} on the \href{http://www.math.sunysb.edu/dynamics/theses}{Stony Brook
  Thesis Server}.

\bibitem[FRS]{markuslassedierk}
Markus F\"orster, Lasse Rempe, and Dierk Schleicher, \emph{Classification of
  escaping exponential maps}, Preprint, 2003,
  \mbox{\href{http://www.arXiv.org/abs/math.DS/0311427}{arXiv:math.DS/0311427}%
}.

\bibitem[FS]{markusdierk}
Markus F\"orster and Dierk Schleicher, \emph{Parameter rays for the exponential
  family}, Preprint, 2005,
  \mbox{\href{http://www.arXiv.org/abs/math.DS/0505097}{arXiv:math.DS/0505097}%
}.


\bibitem[MSS]{mss}
Ricardo Ma{\~n}{\'e}, Paulo Sad, and Dennis Sullivan, \emph{On the dynamics of
  rational maps}, Ann. Sci. \'Ecole Norm. Sup. (4) \textbf{16} (1983), no.~2,
  193--217.

\bibitem[M1]{jackdynamics}
John Milnor, \emph{Dynamics in one complex variable}, Friedr. Vieweg \& Sohn,
  Braunschweig, 1999,
  \mbox{\href{http://www.arXiv.org/abs/math.DS/9201272}{arXiv:math.DS/9201272}%
}.

\bibitem[M2]{jackrays}
\bysame, \emph{Periodic orbits, externals rays and the {M}andelbrot set: an
  expository account}, Ast\'erisque (2000), no.~261, xiii, 277--333,
  G\'eom\'etrie complexe et syst\`emes dynamiques (Orsay, 1995),
 \mbox{\href{http://www.arXiv.org/abs/math.DS/9905169}{arXiv:math.DS/9905169}%
}.



\bibitem[R1]{thesis}
\bysame, \emph{Dynamics of exponential maps}, doctoral thesis,
  Christian-Albrechts-Universit\"at Kiel, 2003,
  \href{http://e-diss.uni-kiel.de/diss_781}{{\tt
  http://e-diss.uni-kiel.de/diss\_781/}}; to appear in the Stony Brook IMS
  preprint series.

\bibitem[R2]{topescapingnew}
\bysame, \emph{Topological dynamics of exponential maps on their escaping
  sets}, Preprint, 2003,
  \mbox{\href{http://www.arXiv.org/abs/math.DS/0309107}{arXiv:math.DS/0309107}%
}, submitted for publication.

\bibitem[R3]{landingsiegel}
\bysame, \emph{{S}iegel disks and periodic rays of entire functions}, Preprint,
  2004,
  \mbox{\href{http://www.arXiv.org/abs/math.DS/0408041}{arXiv:math.DS/0408041}%
}, submitted for publication.

\bibitem[R4]{nonlanding}
\bysame, \emph{Nonlanding dynamic rays of exponential maps}, Preprint, 2005,
  \mbox{\href{http://www.arXiv.org/abs/math.DS/0511588}{arXiv:math.DS/0511588}%
}.


\bibitem[R5]{accessiblesingularvalue}
 \bysame, \emph{On entire functions with accessible singular value}, in
  preparation.


\bibitem[RS1]{boundary}
Lasse Rempe and Dierk Schleicher, \emph{Bifurcations in the space of
  exponential maps}, Preprint \#2004/03, Institute for Mathematical Sciences,
  SUNY Stony Brook, 2004,
  \mbox{\href{http://www.arXiv.org/abs/math.DS/0311480}{arXiv:math.DS/0311480}%
}.

\bibitem[RS2]{expcombinatorics}
\bysame, \emph{Combinatorics of bifurcations in exponential parameter space},
  Preprint, 2004,
  \mbox{\href{http://www.arXiv.org/abs/math.DS/0408011}{arXiv:math.DS/0408011}%
}, submitted for publication.


\bibitem[S1]{habil}
Dierk Schleicher, \emph{On the dynamics of iterated exponential maps},
  Habilitation thesis, TU M\"unchen, 1999.

\bibitem[S2]{expattracting}
\bysame,
  \emph{\href{http://www.math.helsinki.fi/Annales/Vol28/schleich.html}{Attract%
ing dynamics of exponential maps}}, Ann. Acad. Sci. Fenn. Math. \textbf{28}
  (2003), 3--34.

\bibitem[SZ1]{expescaping}
Dierk Schleicher and Johannes Zimmer, \emph{Escaping points of exponential
  maps}, J. London Math. Soc. (2) \textbf{67} (2003), no.~2, 380--400.

\bibitem[SZ2]{expper}
\bysame,
  \emph{\href{http://www.math.helsinki.fi/Annales/Vol28/zimmer.html}{Periodic
  points and dynamic rays of exponential maps}}, Ann. Acad. Sci. Fenn. Math.
  \textbf{28} (2003), 327--354.

\end{thebibliography}

\end{document}